\ifx\shlhetal\undefinedcontrolsequence\let\shlhetal\relax\fi

\documentstyle[12pt]{article}

\input mssymb

\font\teneuf=eufm10
\font\seveneuf=eufm7
\font\fiveeuf=eufm5

\newfam\frakturfam
\textfont\frakturfam\teneuf
\scriptfont\frakturfam\seveneuf
\scriptscriptfont\frakturfam\fiveeuf

\def\frak{\fam\frakturfam\teneuf}

\newcommand{\forces}{\Vdash} 
\newcommand{\decides}{\mathrel{{\Vert}}}

\newcommand{\V}{{\bf V}} 

\newcommand{\lesdot}{\mathrel{\mathord{<}\!\!\raise 
0.8 pt\hbox{$\scriptstyle\circ$}}}

\newcommand{\Q}{{\Bbb Q}} 
\newcommand{\p}{{\Bbb P}}
\newcommand{\B}{{\cal B}}  
\newcommand{\add}{{\rm {\bf add}}\/}

\newcommand{\con}{{\frak c}} 
\renewcommand{\d}{{\frak d}} 
\renewcommand{\b}{{\frak b}}

\newcommand{\cf}{{\rm cf}\/} 
\newcommand{\can}{2^{\textstyle \omega}} 
\newcommand{\fs}{2^{\textstyle <\!\omega}} 
\newcommand{\baire}{\omega^{\textstyle \omega}} 
 
\newcommand{\is}{[\omega]^{\textstyle \omega}} 

\newcommand{\lh}{{\rm lh}\/}

\newcommand{\rest}{{\mathord{\restriction}}} 

\renewcommand{\root}{{\rm root}\/} 
\newcommand{\suc}{{\rm succ}} 

\newcommand{\dom}{{\rm dom}} 
\newcommand{\rng}{{\rm rng}}

\newcommand{\IS}{{\rm IS}\/} 
\newcommand{\CR}{{\rm CR}\/}
\newcommand{\A}{{\cal A}}  
  
\newcommand{\F}{{\cal F}} 
 
\renewcommand{\H}{{\cal H}}  
  
\renewcommand{\L}{{\cal L}}

\renewcommand{\S}{{\Bbb S}}
\newcommand{\s}{{\cal S}}
\newcommand{\T}{{\cal T}} 

\newcommand{\X}{{\cal X}}   
\newcommand{\QED}{\hfill\hspace{0.2in}\vrule width 6pt height 6pt depth 0pt
\vspace{0.1in}} 

\newcommand{\Proof}{\noindent{\sc Proof:} \hspace{0.2in}}
\newtheorem{theorem}{Theorem}[section] 
\newtheorem{lemma}[theorem]{Lemma} 
 
\newtheorem{corollary}[theorem]{Corollary} 
 
\newtheorem{definition}[theorem]{Definition}
\newtheorem{claim}{Claim}[theorem]

\title{More forcing notions imply diamond} 
\author{{\bf Andrzej Ros\l anowski}\thanks{Research partially supported by KBN
654/2/91}\\ 
Dept. of Mathematics and Computer Science\\ 
Bar-Ilan University\\ 
52900 Ramat-Gan, Israel\\ 
and\\ 
Mathematical Institute of Wroclaw University\\ 
50384 Wroclaw, Poland
\and
{\bf Saharon Shelah}\thanks{ Research partially supported by ``Basic
Research Foundation'' administered by The Israel Academy of Sciences and
Humanities. Publication 475.}\\
Institute of Mathematics\\
The Hebrew University of Jerusalem\\
Jerusalem, Israel\\
and\\
Department of Mathematics\\
Rutgers University\\
New Brunswick, NJ, USA}
 
\date\today 

\begin{document} 
\maketitle 
\begin{abstract}
We prove that the Sacks forcing collapses the continuum onto $\d$, answering
the question of Carlson and Laver. Next we prove that if a proper forcing of
the size at most continuum collapses $\omega_2$ then it forces
$\diamondsuit_{\omega_{1}}$.
\end{abstract}

\setcounter{section}{-1}
\section{Introduction} 
In 1979 Baumgartner and Laver proved that after adding $\omega_2$ Sacks
reals (by the countable support iteration) to a model of CH one gets a
model in which the Sacks forcing forces CH (see theorem 5.2 of \cite{BL}).
The question arose when the Sacks forcing may collapse cardinals and which
of them. In 1989 Carlson and Laver posed a hypothesis that the Sacks forcing
collapses the continuum at least onto the dominating number $\d$ (see
\cite{CL}). In the same paper they proved that, assuming CH, the Sacks
forcing forces $\diamondsuit_{\omega_1}$. 
In the present paper we give an affirmative answer to the question of
Carlson and Laver proving that the continuum is collapsed at least onto a
cardinal number called here $\b^{+\epsilon}$ when a Sacks real is added. The
cardinal $\b^{+\epsilon}$ is one of the cardinal invariants laying between
the unbounded number $\b$ and the dominating number $\d$ which were
introduced in \cite{RoS}. After we got the answer we proved that if
$\b^{+\epsilon}=\omega_1$ then the Sacks forcing forces
$\diamondsuit_{\omega_1}$. That naturally suggested the question if this is
an accident and the answer we obtained says that it is a reflection of a more
general theorem. 

The main result of this paper says that if a proper forcing notion $\p$ of
size not greater than the continuum collapses $\omega_2$ then
$\forces_{\p}\diamondsuit_{\omega_1}$. 
\medskip

\noindent{\bf Notation:}\hspace{0.15in}
Our notation is rather standard and is compatible with that of
\cite{Je} or \cite{Ku}. However, there are some exceptions. In a forcing
notion $\p$ we write $p\leq q$ to say that ``the condition $q$ is stronger
than $p$''. The canonical $\p$-name for a generic filter is denoted by
$\Gamma_\p$ or just $\Gamma$.  For a formula $\varphi$ of the forcing
language and a condition $p\in\p$ we say that $p$ {\em decides} $\varphi$
($p\decides\varphi$) if either $p\forces\varphi$ or $p\forces\neg\varphi$.

A forcing notion $(\p,\leq)$ satisfies the Axiom A of Baumgartner (see
\cite{Ba}) if there are partial orders $\leq_n$ on $\p$ (for $n\in\omega$)
such that
\begin{enumerate}
\item $p\leq_0 q$ if and only if $p\leq q$
\item if $p\leq_{n+1}q$ then $p\leq_n q$
\item if a sequence $\langle p_n: n\in\omega\rangle\subseteq\p$ satisfies
$(\forall n\in\omega)(p_n\leq_n p_{n+1})$ then there exists a condition
$p\in\p$ such that $(\forall  n\in\omega)(p_n\leq p)$.
\item if $\A\subseteq\p$ is an antichain, $p\in\p$, $n\in\omega$ then there
exists a condition $q\in\p$ such that $p\leq_n q$ and the set $\{r\in\A:
q\mbox{ and } r \mbox{ are compatible}\}$ is countable.
\end{enumerate}
It is well known that if $\p$ satisfies the Axiom A then $\p$ is proper.

The size of the continuum is denoted by $\con$. We will use the quantifiers
$(\forall^\infty n)$ and $(\exists^\infty n)$ as abbreviations for 
\[(\exists m\in\omega)(\forall n>m)\ \mbox{ and }\ (\forall
m\in\omega)(\exists n>m),\] 
respectively. The Baire space $\baire$ of all functions from $\omega$ to
$\omega$ is endowed with the partial order $\leq^*$:
\[f\leq^* g\ \iff\ (\forall^\infty n)(f(n)\leq g(n)).\]
A family $F\subseteq\baire$ is unbounded in $(\baire,\leq^*)$ if
\[\neg(\exists g\in\baire)(\forall f\in F)(f\leq^* g)\]
and it is dominating in $(\baire,\leq^*)$ if
\[(\forall g\in\baire)(\exists f\in F)(g\leq^* f).\]
The unbounded number $\b$ is the minimal size of an unbounded
family in the partial order $(\baire,\leq^*)$, the dominating number
$\d$ is the minimal size of a dominating family in that order (for more
information about these cardinals see \cite{Va} or \cite{RoS}).

The set of all infinite subsets of $\omega$ is denoted by $\is$. A tree on
$\X$ is a set of finite sequences $T\subseteq \X^{\textstyle <\omega}$ such
that $s\subseteq t\in T$ implies $s\in T$. A tree $T$ on $\X$ is perfect if
for each $s\in T$ there are $t_0,t_1\in T$ extending $s$, both in $T$ and
such that neither $t_0\subseteq t_1$ nor $t_1\subseteq t_0$. The body $[T]$
of a tree $T$ is the set $\{x\in\X^{\textstyle\omega}: (\forall
n\in\omega)(x\rest n\in T)\}$.
\bigskip

\noindent{\bf Acknowledgements:}\hspace{0.15in} Special thanks are due
to the referee for very valuable comments.
\medskip

\section{Antichains of skew trees}
The Sacks forcing $\S$ consists of all perfect trees $T\subseteq\fs$. These
trees are ordered by inclusion (a stronger tree is the smaller one). For
$T\in\S$ and $t\in T$ we say that $t$ {\em ramifies in} $T$ (or $t$ is {\em
a ramification point in} $T$) whenever both $t\hat{\ }0$ and $t\hat{\ }1$
are in $T$. For $s\in T\cap 2^{\textstyle n}$, $n<k$ we say that $s$ {\em
ramifies in $T$ below} $k$ if there is $t\in T$ of length less than $k-1$
such that $s\subseteq t$ and $t$ ramifies in $T$. A node $t\in T$ is a
ramification point of rank $n$ in $T$ if $t$ ramifies in $T$ and exactly $n$
initial segments of $t$ ramify in $T$. Orders $\leq_n$ on $\S$ are defined
by
\begin{quotation}
\noindent $T\leq_n T'$ if and only if 

\noindent $T\leq T'$ and if $t\in T$ is a ramification point of the
rank $\leq n$ then $t\in T'$. 
\end{quotation}
The Sacks forcing $\S$ together with orders $\leq_n$ (for $n\in\omega$)
satisfies Axiom A of Baumgartner (see \cite{Ba}).
\medskip

For $T\in\S$ and $t\in T$ we put $(T)_t=\{s\in T: s\subseteq t\mbox{ or
} t\subseteq s\}$.

\begin{definition}
A tree $T\in\S$ is {\em skew} if for each $n\in \omega$ at most one node
from $T\cap 2^{\textstyle n}$ ramifies in T.
\end{definition}
Clearly the set of all skew perfect trees is dense in $\S$.
\bigskip

Carlson and Laver proved that CH implies
$\forces_{\S}\diamondsuit_{\omega_{1}}$. A detailed analysis of their proof
shows that the result can be formulated as follows.
\begin{theorem}
[T.Carlson, R.Laver, \cite{CL}]
\label{thcl}
\hspace{0.15in} Assume that $\b=\omega_1$ and every maximal antichain
$\A\subseteq\S$ consisting of skew trees is of the size $\con$.\\
Then $\forces_{\S}\diamondsuit_{\omega_{1}}$. \QED
\end{theorem}
\medskip

\noindent Since skew trees are very small (e.g. their bodies are both
meager and null) the question appeared if the second assumption is
always satisfied. The answer is negative:
\begin{theorem}
\hspace{0.15in}It is consistent that there exists a maximal antichain
$\{T_\alpha:\alpha<\omega_{1}\}\subseteq\S$ such that each tree $T_\alpha$
is skew while $\omega_1<\con$.
\end{theorem}

\Proof Let $\bar{T}=\langle T_\alpha:\alpha<\alpha_0\rangle$ be a sequence of
skew trees, $\alpha_0<\omega_{1}$ and let $\s=\{T\in\S:(\forall
\alpha<\alpha_0)(T_\alpha,T\mbox{ are incompatible })\}$. We define a
forcing notion $\Q(\bar{T})$:
\smallskip

\noindent {\em Conditions} are triples $(n,F,\bar{S})$ such that 
\begin{quotation}
\noindent $F\subseteq 2^{\textstyle \leq\! n}$ is a finite skew tree of
height $n\in \omega$, 

\noindent $\bar{S}=\langle S_t: t\in F\cap 2^{\textstyle n}\rangle$,
$t\subseteq\root(S_t)$ and $S_t\in\s$.
\end{quotation}

\noindent {\em The order} is defined by
\begin{quotation}
\noindent $(n^0,F^0,\bar{S}^0)\leq(n^1,F^1,\bar{S}^1)$ if and only if

\noindent $F^1\rest n^0 = F^0$ and $(\forall t\in F^0\!\cap\! 2^{\textstyle
n_0})(\exists s\in F^1\!\cap\! 2^{\textstyle n_1})(S^1_s=(S^0_t)_s)$
\end{quotation}

\begin{claim}
\label{cl1}
\hspace{0.15in}The forcing notion $\Q(\bar{T})$ satisfies the ccc.
\end{claim}
\noindent Why? Suppose that $\langle (n^i,F^i,\bar{S}^i):i<
\omega_{1}\rangle\subseteq \Q(\bar{T})$. First we find $A\in
[\omega_{1}]^{\textstyle \omega_{1}}$, $n\in \omega$ and a finite skew tree
$F\subseteq 2^{\textstyle\leq\! n}$ of the height $n$ such that for each
$i\in A$ we have $n^i=n, F^i=F$. Next we find $A'\in[A]^{\textstyle
\omega_{1}}$, $n^*>n$ and a finite skew tree $F^*\subseteq 2^{\textstyle
n^*}$ such that $F^*\cap 2^{\textstyle n}=F\cap 2^{\textstyle n}$ and for
each $i\in A'$ 
\begin{quotation}
each node $t\in F\cap 2^{\textstyle n}$ ramifies in $F^*$ (below $n^*$) and 

$S^i_t\cap 2^{\textstyle \leq\! n^*}\supseteq (F^*)_t$. 
\end{quotation}
For each $t\in F\cap 2^{\textstyle n}$ choose two distinct $l(t), r(t)\in
(F^*)_t\cap 2^{\textstyle n^*}$. Let $i,j\in A'$. For $t\in F\cap
2^{\textstyle n}$ put $S^*_{l(t)}=(S^i_t)_{l(t)}$ and
$S^*_{r(t)}=(S^j_t)_{r(t)}$. Clearly $(n^*,F^*,\bar{S}^*)\in\Q(\bar{T})$ and
this condition is stronger than both $(n^i,F^i,\bar{S}^i)$ and
$(n^j,F^j,\bar{S}^j)$. The claim is proved. 
\medskip

Suppose that $G\subseteq\Q(\bar{T})$ is a generic filter over $\V$. Then
a density argument shows that $T_G=\bigcup\{F:(\exists
n,\bar{S})((n,F,\bar{S})\in G)\}$ is a skew perfect tree. Let
$\dot{T}_\Gamma$ be the canonical $\Q(\bar{T})$-name for the tree $T_G$. 

\begin{claim}
\label{cl2}
\hspace{0.15in} If $(n,F,\bar{S})\in\Q(\bar{T})$,
$t\in F\cap 2^{\textstyle n}$\\
then $(n,F,\bar{S})\forces \mbox{``}\dot{T}_\Gamma, S_t$ are compatible''.
\end{claim}
\noindent Why? Suppose $(n^0,F^0,\bar{S}^0)\in\Q(\bar{T})$, $t_0\in F^0\cap
2^{\textstyle n^0}$. Take $n^1$ such that $t_0$ ramifies in $S^0_{t_0}$
below $n^1$. Take two distinct extensions $t_0^0,t_0^1$ of $t_0$,
$t^0_0,t^1_0\in S^0_{t_0}\cap 2^{\textstyle n^1}$ and for $t\in(F^0\cap
2^{\textstyle n^0})$ fix an extension $t^1\supseteq t$, $t^1\in S^0_t$. Put
\[F^1=\{t^1\rest m:m\leq n^1\}\cup\{t^i_0\rest m: m\leq n^1, i=0,1\},\
S^1_{t^1}=(S^0_t)_{t^1},\  S^1_{t^i_0}=(S^0_{t_0})_{t^i_0}.\]
Then $(n^1, F^1,\bar{S}^1)\in\Q(\bar{T})$ is a condition stronger than
$(n^0,F^0,\bar{S}^0)$ and
\[(n^1,F^1,\bar{S}^1)\forces t^0_0,t^1_0\in \dot{T}_\Gamma\cap S^0_{t_0}.\] 
Since $S^1_{t^0_0},S^1_{t^1_0}\subseteq S^0_{t_0}$ easy density
argument proves the claim.

\begin{claim}
\label{cl3}
\hspace{0.15in}$\forces_{\Q(\bar{T})}(\forall
\alpha<\alpha_0)(T_\alpha,\dot{T}_\Gamma\mbox{ are incompatible})$.
\end{claim}
\noindent Why? Let $\alpha<\alpha_0$, $(n,F,\bar{S})\in\Q(\bar{T})$. Since each
$S_t$ (for $t\in F\cap 2^{\textstyle n}$) is incompatible with $T_\alpha$ we
find $n^*>n$ and $v(t)\in S_t\cap 2^{\textstyle n^*}$ for $t\in F\cap
2^{\textstyle n}$ such that $v(t)\notin T_\alpha$. Let 
\[F^*=\{v(t)\rest m: m\leq n^*, t\in F\cap 2^{\textstyle n}\}\mbox{ and }
S^*_{v(t)}=(S_t)_{v(t)}\mbox{ for } t\in F\cap 2^{\textstyle n}.\]
Then $(n^*,F^*,\bar{S}^*)\geq (n,F,\bar{S})$ and $(n^*,F^*,\bar{S}^*)\forces
\dot{T}_\Gamma\cap T_\alpha\subseteq F^*$. The claim is proved.

Now we start with $\V\models\neg$CH. Let
$\langle\p_\alpha,\dot{\Q}_\alpha:\alpha<\omega_{1}\rangle$ be the finite
support iteration such that
\[\forces_\alpha\dot{\Q}_\alpha=\Q(\langle\dot{T}_\beta:
\beta<\alpha\rangle)\] 
where $\dot{T}_\beta$ is the $\p_{\beta+1}$-name for the generic tree added by
$\dot{\Q}_\beta$. Let $G\subseteq\p_{\omega_{1}}$ be a generic over $\V$.
Since $\p_{\omega_{1}}$ satisfies ccc (by claim \ref{cl1}) we have
$\V[G]\models\neg$CH. By claim \ref{cl3},
$\langle\dot{T}^G_\alpha:\alpha<\omega_{1}\rangle$ is an antichain in $\S$.
We claim that it is a maximal antichain (in $\V[G]$). 

\noindent Suppose that $\dot{T}$ is a $\p_{\omega_{1}}$-name for an element
of $\S$. Then $\dot{T}$ is a $\p_\alpha$-name for some $\alpha<\omega_{1}$.
Assume that $p\in\p_{\omega_{1}}$ is such that 
\[p\forces_{\omega_{1}}(\forall
\alpha<\omega_{1})(\dot{T},\dot{T}_\alpha\mbox{ are incompatible}).\]
Take $\alpha_0>\alpha$ such that $p\in\p_{\alpha_0}$. Since
\[p\forces_{\alpha_0}(\forall
\alpha<\alpha_0)(\dot{T},\dot{T}_\alpha\mbox{ are incompatible})\]
we can extend $p$ to $q=p\cup\{(\alpha_0,(0,\{\emptyset\},
\langle\dot{T}\rangle))\}\in\p_{\omega_{1}}$. It follows from claim
\ref{cl2} that $q\forces_{\omega_{1}}\mbox{``}\dot{T}_{\alpha_0},\dot{T}$ are
compatible'' - a contradiction. The theorem is proved.\QED

\section{When Sacks forcing forces CH}

In this section we show that if $\d=\omega_{1}$ then $\forces_{\S}$CH, and
hence applying the result of the next section we will be able to conclude
$\forces_{\S}\diamondsuit_{\omega_{1}}$ provided $\d=\omega_1$.
[To be more precise, if CH holds then $\forces_{\S}\diamondsuit_{\omega_1}$
by theorem \ref{thcl} of Carlson and Laver. If we are in the situation of
$\d=\omega_1<\con$ then, by corollary \ref{answer}, $\con$ is
collapsed to $\omega_1$ and hence $\omega_2$ is collapsed (by forcing
with $\S$). Now theorem \ref{diamond2} applies.] This answers
the question of T.Carlson and R.Laver (see \cite{CL}).
\medskip

We start with the following general observation.
\begin{lemma}
\label{generallemma}
\hspace{0.15in} Let $\p$ be a forcing notion, $\kappa$ a cardinal. Suppose
that there exist antichains $\A_\zeta\subseteq\p$ for $\zeta<\kappa$ such
that
\begin{description}
\item[{\bf (*)}]\hspace{0.6in}$(\forall p\in\p)(\exists
\zeta<\kappa)(|\{q\in\A_\zeta: p\leq q\}|=|\p|)$. 
\end{description}
Then $\forces_{\p}|\p^{\V}|\leq\kappa$.
\end{lemma}

\Proof For each $\zeta<\kappa$, by an easy induction, one can construct
a function $\phi_\zeta:\A_\zeta\longrightarrow\p$ such that for every
$p,p'\in\p$
\begin{quotation}
if $|\{q\in\A_\zeta: p\leq q\}|=|\p|$

then $\phi_\zeta(q)=p'$ for some $q\in\A_\zeta$, $q\geq p$.
\end{quotation}
Now let $\dot{\phi}$ be a $\p$-name for a function from $\kappa$ into
$\p^{\V}$ such that
\[q\forces\dot{\phi}(\zeta)=\phi_\zeta(q)\ \ \ \mbox{ for }\zeta<\kappa,\ 
q\in\A_\zeta.\]
Clearly for each $p,p'\in\p$, if
$\zeta<\kappa$ witnesses {\bf (*)} for $p$ then there is $q\geq p$ such that
$q\forces\dot{\phi}(\zeta)=p'$. Consequently
$\forces_{\p}\rng(\dot{\phi})=\p^{\V}$ and we are done. \QED

Thus to prove that the Sacks forcing collapses continuum we will construct
the respective sequence of antichains in $\S$. The sequence will be produced
from a special family of subsets of $\is$

For a set $X\in\is$ let $\mu_X:\omega\stackrel{\rm onto}{\longrightarrow} X$
be the increasing enumeration of the set $X$. 
\begin{definition}
\begin{description}
\item[(1) ] A family $\F\subseteq\is$ is {\em dominating in} $\is$ if
\[(\forall Y\!\in\!\is)(\exists X\!\in\!\F)(\forall^\infty
n)(|[\mu_X(n),\mu_X(n+1)) \cap Y|\geq 2).\] 
\item[(2) ] A family $\F\subseteq\is$ is {\em weakly dominating in } $\is$
if for every set $Y\in\is$
\[(\exists X\!\in\!\F)(\exists^\infty i)(\forall
j\!<\!2^i)(|[\mu_X(2^i\!+\!j),\mu_X(2^i\!+\!j\!+\!1))\cap Y|\geq 2).\]
\item[(3) ] $\b^{+\epsilon}=\min\{|\F|:\F\subseteq\is\mbox{ is weakly
dominating}\}$. 
\end{description}
\end{definition}
\medskip

\noindent{\sc Remarks:} \hspace{0.15in} 1) Note that if $\F$ is a
dominating family in $\is$ then $\{\mu_X: X\in\F\}$ is a dominating family in
the order $\langle\baire, \leq^*\rangle$. And conversely, if $F\subseteq\baire$
is a dominating family of increasing functions,
$X_{f,n}=\{f(n),f(f(n)),f(f(f(n)))\ldots\}$ (for
$f\!\in\! F, n\!\in\! \omega$) then $\{X_{f,n}\!:f\!\in\! F, n\!\in\!
\omega\}$ is a dominating family in $\is$. In particular the minimal size of
a dominating family in $\is$ is the dominating number $\d$. Clearly each
dominating family is weakly dominating. Consequently $\b^{+\epsilon}\leq\d$.

\noindent 2) We have the following inequalities:
\[\b\leq\b^{+\epsilon}\leq\min\{|X|: X\subseteq\can
\mbox{ is not meager}\}.\]
Moreover, the inequality $\b<\b^{+\epsilon}$ is consistent with ZFC (see
\cite{RoS}; $\b^{+\epsilon}$ is the cardinal
$\d(S_{+\epsilon})$ of that paper).

\noindent 3) One can replace ``$\geq 2$'' in the definition
of a weakly dominating family (and $\b^{+\epsilon}$) by ``$\geq 1$''
(and replace the function $i\mapsto 2^i$ by any other increasing
function) and still the results of this section could be carried on (with this
new $\b^{+\epsilon}$). The reason why we use this definition of
$\b^{+\epsilon}$ is that it fits to a more general schema of cardinal
invariants studied in \cite{RoS}. For example note that the unbounded
number $\b$ equals to
\[\min\{|\F|: \F\subseteq\is\ \&\ (\forall Y\!\in\!\is)
(\exists X\!\in\!\F)(\exists^\infty i)(|[\mu_X(i),\mu_X(i+1))\cap Y|\geq 2)\}\]
and ``$\geq 2$'' in the above cannot be replaced by ``$\geq 1$''.
\medskip

\begin{definition}
Let $T\in\S$, $X\in\is$. We say that the condition $T$ {\em weakly obeys}
the set $X$ if  
\[\hspace{-0.1in}(\exists^\infty i)(\forall j\!<\!2^i)(\forall
t\!\in\! T\cap 2^{\mu_X(2^i+j)})(t \mbox{ ramifies in $T$ below }
\mu_X(2^i+j+1)).\]   
\end{definition}

\begin{lemma}
\label{antichain}
\hspace{0.15in}Suppose $X\in\is$. Then there exists an
antichain $\A\subseteq\S$ such that  

\noindent $({\bf *}_{X})$\hspace{0.7in}if $T\in\S$ weakly obeys $X$ then
$|\{S\in\A: T\leq S\}|=\con$. 
\end{lemma}

\Proof Let $\{T_\alpha:\alpha<\con\}=\{T\in\S: T\mbox{ weakly obeys
}X\}$ be an enumeration with $\con$ repetitions. Let $\{h_\alpha:
\alpha<\con\}\subseteq\baire$ be a family of functions such that
\[(\forall \alpha\!<\!\con)(\forall i\!<\!\omega)(h_\alpha(i)<2^i)\mbox{
and } (\forall \alpha\!<\!\beta\!<\!\con)(\forall^\infty i)(h_\alpha(i)\neq
h_\beta(i)).\] 
Since $T_\alpha$ weakly obeys $X$ we have that for infinitely many $i$, for
each $j<2^i$ each node $t\in T_\alpha\cap 2^{\mu_X(2^i+j)}$ ramifies in
$T_\alpha$ below $\mu_X(2^i+j+1)$. Consequently, for each $\alpha<\con$ we
can construct a condition $S_\alpha\geq T_\alpha$ such that for every
$i\in\omega$:
\begin{quotation}
\noindent if some $t\in S_\alpha\cap 2^{\mu_X(2^i+j)}$ ramifies below 
$\mu_X(2^i+j+1)$, $j<2^i$

\noindent then $j=h_\alpha(i)$. 
\end{quotation}
Note that $(\forall^\infty n)(h_\alpha(n)\neq h_\beta(n))$ implies that
conditions $S_\alpha,S_\beta$ are incompatible. Thus $\A=\{S_\alpha:
\alpha<\con\}$ is an antichain. Clearly this $\A$ works. \QED  

\begin{theorem}
\label{d}
\hspace{0.15in} $\forces_{\S} \con=|(\b^{+\epsilon})^{\V}| $.
\end{theorem}

\Proof Since $\forces_{\S} \con=|\con^{\V}| $ it is enough to show that 
\[\forces_{\S}\mbox{``there exists a function }\phi\mbox{ from
}(\b^{+\epsilon})^{\V} \mbox{ onto }\con^{\V}\mbox{''}.\] 
By the definition of the cardinal $\b^{+\epsilon}$ there exists a
sequence $\langle X_\zeta:\zeta<\b^{+\epsilon}\rangle\subseteq\is$ which
is weakly dominating. Apply lemma~\ref{antichain} to
construct antichains $\A_\zeta\subseteq\S$ such that
 
if $T\in\S$ weakly obeys $X_\zeta$ then $|\{S\in\A_\zeta: T\leq S\}|=\con$.

\noindent Since each tree $T\in\S$ weakly obeys some $X_\zeta$ we can
conclude the assertion from lemma~\ref{generallemma}. \QED 

\begin{corollary}
\label{answer}
\hspace{0.15in}Assume that $\d=\omega_{1}$. Then $\forces_{\S}$CH.\QED
\end{corollary}

The Marczewski ideal $\s_0$ is a $\sigma$-ideal of subsets of the Cantor
space $\can$. This ideal is connected with the Sacks forcing. It consist of
all sets $A\subseteq\can$ such that
\[(\forall T\in\S)(\exists T'\geq T)([T']\cap A=\emptyset),\]
where $[T']=\{x\in\can:(\forall n\in \omega)(x\rest n\in T')\}$.

Some connections between the Marczewski ideal $\s_0$ and the Sacks forcing
$\S$ were established in \cite{JMS}.

\begin{corollary}
\label{MarczIdeal}
\hspace{0.15in}$\add(\s_0)\leq\b^{+\epsilon}$
\end{corollary}

\Proof The crucial fact for this inequality is the existence of a
sequence $\langle\A_\zeta^*:\zeta<\b^{+\epsilon}\rangle\subseteq\S$ of
maximal antichains in $\S$ such that
\[(\forall T\in \S)(\exists\zeta<\b^{+\epsilon})([T]\setminus\bigcup
\{[S]:S\in\A_\zeta^*\}\neq\emptyset).\]
For this first, as in the proof of theorem~\ref{d}, find antichains
$\A_\zeta\subseteq\S$ for $\zeta<\b^{+\epsilon}$ such that 
\[(\forall T\in\S)(\exists \zeta<\b^{+\epsilon})(|\{S\in\A_\zeta:T\leq
S\}|=\con).\] 
Now fix $\zeta<\b^{+\epsilon}$. To construct $\A_\zeta^*$ take an enumeration
$\{T_\alpha:\alpha<\con\}$ of $\S$ and an enumeration $\{T_\alpha^*:\alpha<
\con\}$ of $\{T\in\S:|\{S\in\A_\zeta: T\leq S\}|=\con\}$. Next by
induction on $\alpha<\con$ choose trees $S_\alpha\in\S$ and branches
$x_\alpha\in\can$ such that (for $\alpha<\con$):
\begin{quote}
\noindent $x_\alpha\in [T^*_\alpha]\setminus\bigcup_{\beta<\alpha}[S_\beta]$,

\noindent either $(\exists S\in\A_\zeta)(S_\alpha\geq S)$ or $\A_\zeta\cup
\{S_\alpha\}$ is an antichain,

\noindent if $T_\alpha$ is incompatible with all $S_\beta$ (for
$\beta<\alpha$) then $S_\alpha\geq T_\alpha$,

\noindent $S_\alpha$ is incompatible with each $S_\beta$ for $\beta<\alpha$ and

\noindent $[S_\alpha]\cap\{x_\beta:\beta\leq\alpha\}=\emptyset$.
\end{quote}
At stage $\alpha<\con$ we easily find a suitable $x_\alpha\in [T^*_\alpha]$
since continuum many members of $\A_\zeta$ is stronger than $T^*_\alpha$ and
each $S_\beta$ (for $\beta<\alpha$) is either stronger than some member of
$\A_\zeta$ or incompatible with all elements of $\A_\zeta$. (Remember that two
conditions $S,T\in\S$ are incompatible in $\S$ if and only if $[S]\cap
[T]$ is countable.) If the condition $T_\alpha$ is compatible with some
$S_\beta$ for $\beta<\alpha$ then we put $S_\alpha=S_\beta$. Otherwise we
choose $S\in\S$ such that $T_\alpha$ and $S$ are compatible and either
$S\in\A_\zeta$ or $S$ is incompatible with all members of $\A_\zeta$.
As each perfect set contains continuum many disjoint perfect sets we can find
a tree $S_\alpha\geq T_\alpha, S$ such that
$[S_\alpha]\cap\{x_\beta:\beta\leq\alpha\}=\emptyset$. 

\noindent Then $\{S_\alpha:\alpha<\con\}=\A_\zeta^*$ is a maximal antichain
(note that there could be repetitions in $\{S_\alpha:\alpha<\con\}$). The
points $x_\alpha$ (for $\alpha<\con$) witness that no $[T_\alpha^*]$ is
covered by $\bigcup\{[S]: S\in\A_\zeta^*\}$.\medskip

\noindent Now, having antichains $\A^*_\zeta$ as above, we put
$A_\zeta=\bigcup\{[T]:T\in\A_\zeta^*\}$. Since $\A_\zeta^*$ is a
maximal antichain the complement of $A_\zeta$ is in the ideal $\s_0$.
Moreover, for each $T\in\S$ there is $\zeta<\b^{+\epsilon}$ with
$[T]\backslash A_\zeta\neq\emptyset$. Hence
$\bigcup\limits_{\zeta<\b^{+\epsilon}}(\can\backslash A_\zeta)\notin\s_0$. \QED
\bigskip

\noindent{\sc Remark:}\hspace{0.15in} Recently P.~Simon has proved that
in the results of these section one can replace $\b^{+\epsilon}$ by the
unbounded number $\b$.

\section{Collapse $\omega_2$ -- the continuum will fall down}

In this section we will prove that if the Sacks forcing (or any proper
forcing of size $\leq\con$) collapses $\omega_2$ then it forces
$\diamondsuit_{\omega_1}$. First we will give combinatorial tools needed
for the proof. Let us start with fixing some notation.
\medskip

\noindent For an ordinal $\kappa$ by $\IS(\kappa)$ we will denote the set of
finite incresing sequences with values in $\kappa$. $\chi$ stands for a
``sufficiently large'' cardinal, $\H(\chi)$ is the family of all sets
hereditarily of the cardinality less than $\chi$.
\smallskip

\noindent For $\zeta<\omega_1$ let $\zeta=\{e^\zeta_n:n\in \omega\}$ be an
enumeration.
\smallskip

\noindent Let $S^2_i=\{\delta<\omega_2:\cf(\delta)=\omega_i\}$ for $i=0,1$.

\begin{lemma}
[S.Shelah, see 2.3 of \cite{Sh}]
\label{clubguessing}
\hspace{0.15in} There exists a (``club--guessing'') sequence $\bar{C}=\langle
C_\delta:\delta\in S^2_0\rangle$ such that 
\begin{enumerate}
\item $C_\delta\subseteq\delta$, $\sup C_\delta=\delta$,
\item the order type of $C_\delta$ is $\omega$,
\item for every closed unbounded subset $E$ of $\omega_2$ there exist
$\delta\in S^2_0$ such that $C_\delta\subseteq E$.\QED
\end{enumerate}
\end{lemma}
We fix a club--guessing sequence $\bar{C}=\langle C_\delta:\delta\in
S^2_0\rangle$ as in \ref{clubguessing}. For
$\delta\in S^2_0$ let $C_\delta=\{\alpha^\delta_n:n\in \omega\}$
be the increasing enumeration. 

\begin{definition}\hspace{0.15in}
\label{defcreat}
Let $\delta\in S^2_0$ and let $\zeta<\omega_1$ be limit.
\begin{enumerate}
\item A sequence $\langle N_\eta:\eta\in\IS(\omega)\rangle$ is {\em a
semi-$(\zeta,\delta)$-creature} (for the sequence $\bar{C}$) if
\begin{description}
\item[$\alpha)$] $N_\eta$ is a countable elementary submodel of $\H(\chi)$,
$N_\eta\cap \omega_{1}\subseteq\zeta$ and $\bigcup_{n\in\omega}N_{h\rest
n}\cap\omega_1=\zeta$ for every increasing function $h\in\baire$,
\item[$\beta)$] if $\eta\subseteq\nu$ then $N_\eta\prec N_\nu$,
\item[$\gamma)$]  $N_\eta\cap\omega_2\subseteq\alpha^\delta_0 \cup
\bigcup_{n<\lh(\eta)} [\alpha^\delta_{\eta(n)},\alpha^\delta_{\eta(n)+1})$,
\item[$\delta)$] for each $n<\lh(\eta)$ the intersection $N_\eta\cap
[\alpha^\delta_{\eta(n)}, \alpha^\delta_{\eta(n)+1})$ is non empty.
\end{description}
\item Let $\p$ be a forcing notion, $X\in\H(\chi)$. {\em A
$(\zeta,\delta)$-creature for $\p, X$}
is a system $\{(N_\eta,\tau_\eta,k_\eta):\eta\in\IS(\omega)\}$ such that
\begin{description}
\item[$\alpha)$] the sequence $\langle N_\eta:\eta\in\IS(\omega)\rangle$ is
a semi-$(\zeta,\delta)$-creature and\\
$X,\p,\leq_{\p},\omega_2,\omega_{1},\ldots\in N_\emptyset$,
\item[$\beta)$]  $k_\eta\in\omega$, $\{e^\zeta_k:k<k_\eta\}\subseteq N_\eta$,
and for every increasing function $h\in\baire$ the sequence $\langle k_{h\rest
n}:n\in \omega\rangle$ is unbounded,
\item[$\gamma)$] $\tau_\eta$ is a function such that $\dom(\tau_\eta)\in
[\p\times \omega]^{\textstyle \leq\!\omega}$, for each $k\in \omega$ the set
$\{p\in\p: (p,k)\in\dom(\tau_\eta)\}$ is an antichain in $\p$ and
$\rng(\tau_\eta)\subseteq 2$,
\item[$\delta)$] if $\eta\subseteq\nu$ then $k_\eta\leq k_\nu$ and
$\tau_\eta\subseteq\tau_\nu$.
\end{description}
\item Let $\CR^\zeta_\delta(\p,X)$ be the family of all
$(\zeta,\delta)$-creatures for $\p,X$.
\end{enumerate}
\end{definition}
\medskip

\noindent{\sc Remarks:}\hspace{0.15in} {\em 1.}\ \ \ A $\p$-name for a subset
of $\zeta<\omega_1$ can be thought of as a function $\tau$ such that $\rng\tau
\subseteq 2$ and $\dom\tau\subseteq\p\times\zeta$ has the following property:
\begin{quote}
\noindent for each $\xi\in\zeta$ the set $\{p\in\p: (p,\xi)\in\dom\tau\}$ is
an antichain in $\p$
\end{quote}
(and then for $(p,\xi)\in\tau$: $p\forces\xi\in\tau$ if $\tau(p,\xi)=1$ and
$p\forces\xi\notin\tau$ otherwise). If the forcing notion $\p$ is proper
every such a name can be (above each condition) forced to be equal to a
countable name.
\medskip

\noindent{\em 2.}\ \ \ Thus in a $(\zeta,\delta)$-creature
$\{(N_\eta,\tau_\eta,k_\eta):\eta\in\IS(\omega)\}$ for $\p$ the functions
$\tau_\eta$ can be thought of as approximations of a name for a subset of
$\zeta$. Note that we demand no relations between functions $\tau_\eta$ and
models $N_\eta$. The last are only ``side parameters''. The
parameter will decide above which conditions the name is described by the
functions $\tau$ determined by a branch through the creature.
\bigskip

\begin{lemma}
\label{creatures}
\hspace{0.15in} For every $X\in\H(\chi)$ and a closed unbounded set
$D\subseteq\omega_1$ for some $\zeta\in D$ and $\delta\in S^2_0$ there exists a
semi-$(\zeta,\delta)$-creature $\langle N^*_\eta:\eta\in\IS(\omega)\rangle$
such that $X\in N^*_\emptyset$.
\end{lemma}

\Proof The following special case of theorem 2.2 of \cite{RS} is a main
tool for constructing semi-creatures:

\begin{claim}
[M.Rubin and S.Shelah, \cite{RS}]
\label{rubinshelah}
\hspace{0.15in} Suppose that $\T\subseteq\omega_2^{\textstyle
<\omega}$ is a tree such that for each node $t\in\T$ the set $\suc_{\T}(t)$
of successors of $t$ is of the size $\omega_2$. Assume that
$\phi:\T\longrightarrow\omega_1$. Then there exists a subtree $\T_0$ of $\T$
such that
\[(\forall t\in\T_0)(|\suc_{\T_0}(t)|=\omega_2)\ \mbox{ and }\ \sup\phi
[\T_0]<\omega_1.\]
If additionally $\phi$ is increasing (i.e. $t\subseteq s\in\T$ implies
$\phi(t)\leq\phi(s)$) then we can demand that $\lim_n \phi(x\rest n)$ is
constant for all infinite branches $x\in[\T_0]$.
\end{claim}
\medskip

For $v\in\IS(S^2_1)$ choose $N_v$ such that
\begin{description}
\item[(0) ] $X\in N_\emptyset$;
\item[(1) ] $N_v$ is an elementary countable submodel of $\H(\chi)$;
\item[(2) ] $N_v\cap \omega_{1}\in D$, $[\max(v),\omega_2)\cap
N_v\neq\emptyset$; 
\item[(3) ] if $v\subseteq w$ then $N_v\prec N_w$.
\end{description}
Now we will inductively define a tree $\T\subseteq\IS(S^2_1)$ and ordinals
$\delta_v<\omega_2$ for $v\in\T$ such that:
\medskip

\begin{description}
\item[(4) ] if $v\in\T$ then $|\suc_{\T}(v)|=\omega_2$ and
\item[(5) ] if $v\in\T$, $\phi_v:S^2_1\stackrel{\rm onto}{\longrightarrow}
\suc_{\T}(v)$ is the increasing enumeration of $\suc_{\T}(v)$ then for every
$\alpha\in S^2_1$ and $w\in \T$, $w\supseteq v\hat{\ }\phi_v(\alpha)$
\[N_w\cap \alpha\subseteq \delta_v.\]
\end{description}
\medskip

\noindent To start with we put $\emptyset\in\T$. For each $v\in\IS(S^2_1)$ let
$\rho_v=\sup(N_v\cap v(0))$. Applying claim \ref{rubinshelah} for each
$\alpha\in S^2_1$ we find a tree $\T^{\langle \alpha\rangle}\subseteq
\IS(S^2_1)$ and $\rho^\alpha<\alpha$ such that
\begin{description}
\item[(6) ] $\root(\T^{\langle\alpha\rangle})=\langle \alpha\rangle$;
\item[(7) ] each node extending $\langle \alpha\rangle$ has $\omega_2$
successors in $\T^{\langle \alpha\rangle}$;
\item[(8) ] for each $v\in\T^{\langle \alpha\rangle}$, $\rho_v<\rho^\alpha$.
\end{description}
Applying Fodor's lemma we find $\delta_\emptyset$ and $A_\emptyset$ such
that 
\begin{description}
\item[(9) ] $A_\emptyset\in[S^2_1]^{\textstyle \omega_2}$;
\item[(10) ] $\delta_\emptyset=\rho^\alpha$ for $\alpha\in A_\emptyset$.
\end{description}
We put $A_\emptyset=\suc_{\T}(\emptyset)$ and we decide that $(\T)_{\langle
\alpha\rangle}\subseteq\T^{\langle \alpha\rangle}$ for each $\alpha\in
A_\emptyset$. Note that at this moment we are sure that if
$\langle\phi_\emptyset(\alpha)\rangle\subseteq w$ then $N_w\cap
\alpha\subseteq N_w\cap\phi_\emptyset(\alpha)\subseteq\delta_\emptyset$ for
each $\alpha\in S^2_1$. 

Suppose we have decided that $v\in\T$ and $(\T)_v\subseteq\T^v$.

Let $\phi_v':S^2_1\stackrel{\rm onto}{\longrightarrow}\suc_{\T^v}(v)$ be the 
increasing enumeration. For each $\alpha\in S^2_1$ we apply claim
\ref{rubinshelah} to
find $\rho^\alpha<\alpha$ and a tree $\T^{v\hat{\ }\phi_v'(\alpha)}\subseteq
\T^v$ such that 
\begin{description}
\item[(11) ] $\root(\T^{v\hat{\ }\phi_v'(\alpha)})=v\hat{\
}\phi_v'(\alpha)$; 
\item[(12) ] each node in $\T^{v\hat{\ }\phi_v'(\alpha)}$ extending $v\hat{\
}\phi_v'(\alpha)$ has $\omega_2$ successors in $\T^{v\hat{\
}\phi_v'(\alpha)}$; 
\item[(13) ] for each $w\in\T^{v\hat{\ }\phi_v'(\alpha)}$, $w\supseteq
v\hat{\ }\phi_v'(\alpha)$ we have $\sup(N_w\cap \alpha)<\rho^\alpha$. 
\end{description}
Next we choose $\delta_v$ and $A_v$ such that 
\begin{description}
\item[(14) ] $A_v\in [S^2_1]^{\textstyle \omega_2}$;
\item[(15) ] $\delta_v=\rho^\alpha$ for all $\alpha\in A_v$.
\end{description}
We put $\suc_{\T}(v)=\phi_v'[A_v]$ and we decide that for $\alpha\in A_v$ 
\[(\T)_{v\hat{\ }\phi_v'(\alpha)}\subseteq\T^{v\hat{\ }\phi_v'(\alpha)}.\] 
Note that at this moment we are sure that if $w\in\T$, $v\hat{\
}\phi_v(\alpha)\subseteq w$ then $N_w\cap \alpha\subseteq N_w\cap
\beta\subseteq\delta_v$, where $\phi_v'(\beta)=\phi_v(\alpha)$ (clearly
$\alpha\leq \beta$). This finishes the construction of the tree $\T$
(satisfying {\bf (4), (5)}).
\medskip

For $v\in\T$ let $\zeta_v=N_v\cap \omega_{1}\in D$. We apply claim
\ref{rubinshelah} once
again to find $\zeta<\omega_{1}$ and a tree $\T^*\subseteq\T$ such that each
node in $\T^*$ has $\omega_2$ successors in $\T^*$ and for each
$\omega$-branch $z$ through $\T^*$ we have $\sup\{\zeta_{z\rest n}:n\in
\omega\}=\zeta$. Then $\zeta\in D$. For $v\in\T^*$ let $\psi_v:S^2_1\stackrel
{\rm onto}{\longrightarrow}\suc_{\T^*}(v)$ be the increasing enumeration and
let $\delta^*_v=\sup(N_v\cap \omega_2)$. Let
\[\begin{array}{rl}
E=\{\delta\!<\!\omega_2: & \delta\mbox{ is limit } \&\ (\forall
v\!\in\!\IS(\delta)\cap\T^*)(\delta_v<\delta\ \& \ \delta^*_v<\delta)\ \&\\
\& & \ (\forall v\!\in\!\T^*\cap\IS(\delta))(\forall\beta\!<\!\delta)(\exists
\gamma\!\in\!S^2_1)(\beta<\gamma\leq\psi_v(\gamma)<\delta)\}.
\end{array}\]
Since $E$ is a closed unbounded subset of $\omega_2$ we find
$\delta\in E$ such that $C_\delta\subseteq E$. 
\medskip

Now we may define the semi-$(\zeta,\delta)$-creature we are looking for
by constructing an embedding $\pi:\IS(\omega)\longrightarrow\T^*$ such that
$\lh(\pi(\eta))=\lh(\eta)$ and choosing corresponding models
$N_{\pi(\eta)}$. This is done by induction on the length of a sequence
$\eta\in\IS(\omega)$:

\noindent Put $\pi(\emptyset)=\emptyset$. Note that
$\delta_{\emptyset}<\alpha^\delta_0$ (as $\alpha^\delta_0\in E$).

\noindent Suppose we have defined $\pi(\eta)\in\T^*$ such that
$\delta_{\pi(\eta)}<\alpha^\delta_{n_{k-1}+1}$, where $\eta=\langle
n_0,n_1,\ldots,n_{k-1}\rangle$. Given $n_k>n_{k-1}$.

\noindent Take any $\gamma\in
(\alpha^\delta_{n_k},\alpha^\delta_{n_k+1})\cap S^2_1$
and put $\pi(\eta\hat{\ }n_k)=\pi(\eta)\hat{\
}\psi_{\pi(\eta)}(\gamma)\in\T^*$. By the choice of $\gamma$ we have
$\delta_{\pi(\eta\hat{\ }n_k)}<\alpha^\delta_{n_k+1}$. 
\medskip

Finally let $N_\eta^*=N_{\pi(\eta)}$ for $\eta\in\IS(\omega)$.
\medskip

\noindent Since $\delta_\emptyset<\alpha^\delta_0$, $\delta^*_{\pi(\langle
n_0\rangle)}<\alpha^\delta_{n_0+1}$ we have that for every $n_0\in \omega$
\[N^*_{\langle n_0\rangle}\cap \omega_2\subseteq \alpha^\delta_{0}\cup
[\alpha^\delta_{n_0},\alpha^\delta_{n_0+1})\ \mbox{ and }\ 
N^*_{\langle n_0\rangle}\cap
[\alpha^\delta_{n_0},\alpha^\delta_{n_0+1})\neq\emptyset\] 
(we use here (2) and (5)). Similarly, if $\eta=\langle
n_0,\ldots,n_{k-1},n_k\rangle\in\IS(\omega)$ then 
\[N^*_\eta\cap\alpha^\delta_{n_{i+1}}\subseteq\alpha^\delta_{n_i+1}\mbox{ for
} i<k \mbox{ and }\]
\[N^*_\eta\cap[\alpha^\delta_{n_i},\alpha^\delta_{n_i+1})\neq\emptyset\mbox{
for } i\leq k.\] 
Consequently the sequence $\langle N^*_\eta:\eta\in\IS(\omega)\rangle$
is a semi-$(\zeta,\delta)$-creature (and we are done as $X\in N^*_\emptyset$,
$\zeta\in D$). \QED
\bigskip

\begin{theorem}
\label{diamond2}
\hspace{0.15in} Assume $\p$ is a proper forcing notion, $|\p|\leq\con$.
Suppose $\forces_{\p} |\omega_2^{\V}|=\omega_1$. Then
$\forces_{\p} \diamondsuit_{\omega_1}$.
\end{theorem}

\Proof Let $\p$ be a proper forcing notion collapsing $\omega_2$ and of
size $|\p|\leq\con$. Since $\p$ collapses $\omega_2$ and $|\p|\leq\con$
we have $\con\geq\omega_2$. Let $\Theta$ be a $\p$-name such that
\[\forces_{\p}\mbox{``}\Theta: \omega_1\longrightarrow
\omega_2^{\V}\mbox{ is an increasing unbounded function''.}\]

Our aim is to construct a sequence $\langle \dot{A}_\zeta:\zeta<
\omega_1\rangle$ of $\p$-names which witnesses $\diamondsuit_{\omega_1}$
in $\V^{\p}$. In the construction we will use $(\zeta,\delta)$-creatures
which can be thought of as countable ``trees'' of possible
fragments of names for subsets of $\zeta$ (together with some parameters for
controlling their behaviour). Each infinite branch through the creature will
define a (countable) name for a subset of $\zeta$. Next we will choose
continuum many branches together with conditions in $\p$. Our choice will
ensure that the conditions form an antichain in $\p$ and all antichains
involved in the name determined by a single branch (in important cases) are
predense above the corresponding condition. This will define the name
$\dot{A}_\zeta$ for a subset of $\zeta$. The main difficulty will be in
proving that the sequence $\langle\dot{A}_\zeta:\zeta<\omega_1\rangle$ is (a
name for) a $\diamondsuit_{\omega_1}$--sequence. But this we will obtain right
from the existence of creatures which was proved in lemma \ref{creatures}.

Before we define the names $\dot{A}_\zeta$ we have to identify some
creatures (as the set $\CR^\zeta_\delta(\p,\Theta)$ can be very large):
\medskip

\noindent For a $(\zeta,\delta)$-creature $S=\{(N_\eta,\tau_\eta,k_\eta):
\eta\in\IS(\omega)\}\in\CR^\zeta_\delta(\p,\Theta)$ let
\[U(S)=\bigcup_{\eta\in\IS(\omega)} N_\eta\cap\p.\]
Clearly $U(S)$ is a countable subset of $\p$ and hence there is at most $\con$
possibilities for $U(S)$. Let $S^i=\{(N_\eta^i,\tau_\eta^i,k_\eta^i):
\eta\in\IS(\omega)\}\in\CR^\zeta_\delta(\p,\Theta)$, $i=0,1$. We say that the
creatures $S^0,S^1$ are {\em equivalent} ($S^0\equiv S^1$) whenever
\begin{description}
\item[(i) ] $U(S^0)=U(S^1)$ and
\item[(ii) ] for each $\eta\in\IS(\omega)$: $N^0_\eta\cap\p=N^1_\eta\cap\p$,
$k^0_\eta=k^1_\eta$, $\tau^0_\eta=\tau^1_\eta$ and \\
$\{A^0\cap U(S^0):A^0\in N^0_\eta\mbox{ is a maximal antichain in }\p\}=$\\
$=\{A^1\cap U(S^1):A^1\in N^1_\eta\mbox{ is a maximal antichain in }\p\}$.
\end{description}
(Note that actually condition (ii) implies (i).)
Since for each $\eta\in\IS(\omega)$ there is at most $\con$ possibilities
for $k_\eta,\tau_\eta$, $N_\eta\cap\p$ and $\{A\cap U(S):A\in N_\eta$ is a
maximal antichain in $\p\}$  the relation $\equiv$ has at most $\con$
equivalence classes.  

The following claim should be clear:

\begin{claim}
\label{clear}
\hspace{0.15in} Let $S^i=\{(N_\eta^i,\tau_\eta^i,k_\eta^i):
\eta\in\IS(\omega)\}\in\CR^\zeta_\delta(\p,\Theta)$ (for $i=0,1$) be equivalent
creatures. Let $h\in\baire$ be an increasing function. Then
\begin{enumerate}
\item $\bigcup_{n\in\omega}N^i_{h\rest n}$ is an elementary (countable)
submodel of $\H(\chi)$,
\item $\bigcup_{n\in\omega}N^0_{h\rest n}\cap\p=\bigcup_{n\in\omega}
N^1_{h\rest n}\cap\p$,
\item if $A^0\in N^0_{h\rest n}$ is a maximal antichain in $\p$ then
for some maximal antichain $A^1\in N^1_{h\rest n}$ we have $A^0\cap
\bigcup_{n\in\omega}N^0_{h\rest n}=A^1\cap\bigcup_{n\in\omega}N^1_{h\rest n}$,
\item $\{A^0\cap\bigcup_{n\in \omega}N^0_{h\rest n}:A^0\in \bigcup_{n\in
\omega}N^0_{h\rest n}\mbox{ is a maximal antichain in }\p\}=$

$=\{A^1\cap\bigcup_{n\in \omega}N^1_{h\rest n}:A^1\in
\bigcup_{n\in \omega}N^1_{h\rest n}\mbox{ is a maximal antichain in }\p\}$,
\item if $p\in\p$ is $(\bigcup_{n\in \omega}N^0_{h\rest n},\p)$-generic then
it is $(\bigcup_{n\in \omega}N^1_{h\rest n},\p)$-generic.
\end{enumerate}
\end{claim}
\medskip

Fix a limit ordinal $\zeta<\omega_1$.
\medskip

We are going to define a name $\dot{A}_\zeta$ for a subset of $\zeta$.

Suppose that $\bigcup_{\delta\in S^2_0}\CR^\zeta_\delta(\p,\Theta)\neq
\emptyset$.

Let $p^i\in\p$, $S^i\in\bigcup_{\delta\in S^2_0}\CR^\zeta_\delta(\p,
\Theta)$ (for
$i<\con$) be such that $\{(p^i,[S^i]_\equiv):i<\con\}$ lists of all members
of $\p\times (\bigcup_{\delta\in S^2_0}\CR^\zeta_\delta(\p,\Theta)/\!\equiv)$
with $\con$ repetitions. Take any family $\{h_i:i<\con\}\subseteq\baire$ of
increasing functions  such that for distinct $i,j<\con$ the intersection 
$\rng(h_i)\cap\rng(h_j)$ is finite.  

Now for each $i<\con$ we put $M_i=\bigcup_{n<\omega}N^i_{h_i\rest n}$,
$\tau_i=\bigcup_{n\in \omega}\tau^i_{h_i\rest n}$.

Each $M_i$ is a countable elementary submodel of $\H(\chi)$
and $\tau_i$ is a function. Since $\p$ is proper we find $p_i\in\p$ such
that $p_i$ is $(M_i,\p)$-generic. If we can find such a condition $p_i$
above the condition $p^i$ then we also demand $p_i\geq p^i$. 
Note that $M_i\cap \omega_{1}=\zeta$ and
$M_i\cap \omega_2\subseteq\delta^i$ is cofinal in $\delta^i$ (what is a
consequence of $(\gamma)$, $(\delta)$ of definition \ref{defcreat}(1)),
where $\delta^i\in S^2_0$ is such that $S^i\in\CR^\zeta_{\delta_i}(\p,
\Theta)$. Hence
\[p_i\forces\mbox{``}\rng(\Theta\rest\zeta)\subseteq M_i\mbox{ is unbounded
in } \delta^i\mbox{''}.\]

If $\delta^i\neq\delta^j$ then the conditions $p_i,p_j$ force inconsistent
sentences (unboundness of $\rng(\Theta\rest\zeta)$ in $\delta^i$, $\delta^j$,
respectively). If $\delta^i=\delta^j$ but $i\neq j$ then the choice of the
functions $h_i,h_j$ guaranties (by $(\gamma)$, $(\delta)$ of
\ref{defcreat}(1)) that sets $M_i\cap
[\alpha,\delta^i)$ and $M_j\cap [\alpha,\delta^i)$ are disjoint for some
$\alpha<\delta^i$. Consequently if $i\neq j$ then $p_i,p_j$ are incompatible.

Let $\A_\zeta$ be a maximal antichain in $\p$ extending
$\{p_i:i<\con\}$ and let $\dot{A}_\zeta$ be a name for a subset of
$\zeta$ such that for each $(p,k)\in\dom(\tau_i)$
\[p_i\forces\mbox{``if }p\in\Gamma_{\p}\mbox{ then }\dot{A}_\zeta(e^\zeta_k) 
=\tau_i(p,k)\mbox{''}\]  
(we identify a subset of $\zeta$ with its characteristic function).

If $\bigcup_{\delta\in S^2_0}\CR^\zeta_\delta(\p,\Theta)=\emptyset$
then take any maximal antichain and a name for a subset of $\zeta$.

\noindent We want to show that the sequence $\langle\dot{A}_\zeta:
\zeta<\omega_1\rangle$ is a (name for a) $\diamondsuit_{\omega_1}$-sequence.
For this suppose that $\dot{A}$ is a $\p$-name for a subset of
$\omega_{1}$, $\dot{D}$ is a $\p$-name for a closed unbounded subset of
$\omega_1$, $p\in\p$. We have to prove:

\begin{claim}
\hspace{0.15in} There exist a limit ordinal $\zeta<\omega_1$ and a condition
$p^*\in\A_\zeta$ such that $p^*\geq p$ and $p^*\forces$``$\dot{A}_\zeta=
\dot{A}\cap\zeta\ \&\ \zeta\in\dot{D}$''.
\end{claim}

\noindent To prove the claim we use lemma \ref{creatures} to find a
semi-$(\zeta,\delta)$-creature $\langle N^*_\eta: \eta\in\IS(\omega)\rangle$
such that $\Theta,\dot{A},\dot{D},\p,p,\ldots\in N^*_\emptyset$. Next:

let $k_\eta=\min\{l: e^\zeta_l\notin N^*_{\eta}\}$.

For each $\eta\in\IS(\omega)$ and $k<k_\eta$ we fix a maximal antichain
$\B^k_\eta$ in $\p$ such that $\B^k_\eta\in N^*_\eta$ and $(\forall
p\in\B^k_\eta)(p\decides e^\zeta_k\in \dot{A})$. Moreover we demand that
$\eta\subseteq\nu$ implies $\B^k_\eta=\B^k_\nu$ (for $k<k_\eta$). Now we
define functions $\tau_\eta$ for $\eta\in\IS(\omega)$ by
\begin{quotation}
$\dom(\tau_\eta)=\bigcup_{k<k_\eta}((\B^k_\eta\cap N^*_\eta)\times\{k\})$,

$\tau_\eta(p,k)=1$ if and only if $p\forces e^\zeta_k\in\dot{A}$.
\end{quotation}
It should be clear that $S=\{(N^*_\eta,\tau_\eta,k_\eta):\eta\in\IS(\omega)\}$
is a $(\zeta,\delta)$-creature for $\p, \Theta$. Thus we find $i<\con$ such
that $S\equiv S^i$ and $p=p^i$ (where $S^i\in\CR^\zeta_\delta(\p,\Theta)$,
$p^i\in\p$ are as in the definition of the antichain $\A_\zeta$ and the name
$\dot{A}_\zeta$). Then the condition $p^*=p_i\in\A_\zeta$ is $(\bigcup_{n\in
\omega}N^*_{h_i\rest n},\p)$-generic, $p^*\geq p^i=p\in N^*_\emptyset$. The
name $\dot{A}_\zeta$ agrees with decissions of $\tau_{h_i\rest n}$ (or
$\B^k_{h_i\rest n}$). By the genericity of $p^*$ we conclude that
$p^*\forces$``$\dot{A}\cap\zeta=\dot{A}_\zeta\ \&\ \zeta\in\dot{D}$''.
Thus the claim is proved.
\medskip

The theorem follows from the claim. \QED

\section{Laver forcing, Miller forcing, Silver forcing...}

Results of the second section can be formulated for other forcing notions.
Without any problems we can prove the respective facts for the Silver
forcing (and generally for forcing notions consisting of compact trees). 

Recall that the Silver forcing notion consists of partial functions $p$ such
that $\dom(p)\subseteq \omega$, $\omega\backslash\dom(p)$ is infinite and
$\rng(p)\subseteq 2$. These functions are ordered by the inclusion.

\begin{theorem}
The Silver forcing notion forces ``$\con=|(\b^{+\epsilon})^{\V}|$''.\QED
\end{theorem}

We have to be more carefull when we work with trees on $\omega$.
Nevertheless even in this case we get the similar result

The Laver forcing $\L$ consists of infinite trees $T\subseteq
\omega^{\textstyle <\! \omega}$ such that for each $t\in T$,
$\root(T)\subseteq t$ we have $|\suc_T(t)|=\omega$.

\begin{definition}
W say that a condition $T\in\L$ {\em weakly obeys} a set $X\in\is$
whenever for each ramification point $t\in T$ 
\[(\exists^\infty i)(\forall j<2^i)(\suc_T(t)\cap
[\mu_X(2^i+j),\mu_X(2^i+j+1))\neq\emptyset).\]
\end{definition}

Fix $T\in\L$. Take $X_0\in\is$ such that for each ramification point
$t\in T$
\[(\forall^\infty i)(\suc_T(t)\cap
[\mu_{X_0}(i),\mu_{X_0}(i+1))\neq\emptyset).\] 
Suppose that $X\in\is$ is such that 
\[(\exists^\infty i)(\forall
j<2^i)(|[\mu_X(2^i+j),\mu_X(2^i+j+1))\cap X_0|\geq 2).\]
Then clearly $T$ weakly obeys $X$. Consequently if $\F\subseteq\is$ is a
weakly dominating family then $T$ weakly obeys some $X\in\F$. 

Suppose now that $T$ weakly obeys $X\in\is$ and
$h:\omega\longrightarrow \omega$ is such that $(\forall
i)(h(i)<2^i)$. Then we can easily construct a condition $T^h\geq
T$ such that
\begin{quotation}
\noindent if $t\in T^h$ is a ramification point in $T^h$, $t\hat{\
}n\in T^h$ and  $j<2^i$, $2^i+j\leq n<2^i+j+1$

\noindent then $h(i)=j$. 
\end{quotation}
Moreover, if $h_0,h_1$ are such that $(\forall^\infty i)(h_0(i)\neq
h_1(i))$ then the respective conditions $T^{h_0}, T^{h_1}$ are
incompatible -- their intersection has no node with infinitely many
immediate successors. Consequently we can repeat the proof of
\ref{antichain} and we get

\begin{theorem}
\hspace{0.15in} $\forces_{\L}\con=|(\b^{+\epsilon})^{\V}|$.\QED
\end{theorem}

The argument above applies for the Miller forcing too. Recall that
this order consists of perfect trees $T\subseteq[\omega]^{\textstyle
<\!\omega}$ such that 
\[(\forall t\in T)(\exists s\in T)(t\subseteq s\ \&\
|\suc_T(s)|=\omega).\] 
Thus we can conclude

\begin{theorem}
\hspace{0.15in}The Miller forcing collapses the continuum onto
$(\b^{+\epsilon})^{\V}$. \QED
\end{theorem}

\shlhetal

\bigskip

\begin{thebibliography}{20}
\bibitem{Ba} Baumgartner J., {\em Iterated forcing} in Surveys in Set
Theory, ed. by A.R.D.Mathias, London Math. Soc., Lecture Notes 87, pp 1--59.  
\bibitem{BL} Baumgartner J., Laver R., {\em Iterated perfect set forcing},
Annals of Mathematical Logic 17(1979), pp 271--288.
\bibitem{CL} Carlson T., Laver R., {\em Sacks reals and Martin's Axiom},
Fundamenta Mathematica, vol 133(1989), pp 161--168.
\bibitem{Ku} Kunnen K., {\em Set Theory (An introduction to Independence
Proofs)}, North--Holland, Amsterdam, 1983. 
\bibitem{Je} Jech T., {\em Set Theory}, Academic Press, New York 1978.
\bibitem{JMS} Judah H., Miller A., Shelah S., {\em Sacks forcing, Laver
forcing and Martin's Axiom}, Archive for Mathematical Logic, vol 31(1992),
pp 145--161. 
\bibitem{RoS} Roslanowski A., Shelah S., {\em Localizations of subsets of
$\omega$}, submitted to Archive for Mathematical Logic.
\bibitem{RS} Rubin M., Shelah S., {\em Combinatorial problems on trees:
partitions, $\Delta$-systems and large free subtrees}, Annals of Pure and
Applied Logic, vol 33(1987), pp 43--81.
\bibitem{Sh} Shelah S., {\em There are Jonsson algebras in many 
inaccessible cardinals}, in {\bf Cardinal Arithmetic}, Oxford University
Press, in press.
\bibitem{Va} Vaughan J.E., {\em Small uncountable cardinals in topology} in
{\bf Open problems in topology}, J.van Mill and G.M.Reed eds., pp 195--218, 
North-Holland, Amsterdam 1990.

\end{thebibliography}
\end{document}